\documentclass[leqno,11pt]{amsart}
\usepackage{amsmath,amscd}
\usepackage{epsfig,graphics}
\usepackage{amssymb,latexsym}
\usepackage{mathrsfs}
\usepackage{eucal}

\newtheorem{thm}{\bf Theorem}[section]
\newtheorem{df}[thm]{\bf Definition}
\newtheorem{prop}[thm]{\bf Proposition}
\newtheorem{cor}[thm]{\bf Corollary}
\newtheorem{lem}[thm]{\bf Lemma}
\newtheorem{rem}[thm]{\bf Remark}
\newtheorem{ex}[thm]{\bf Example}

\newcommand{\bs}{\boldsymbol}
\newcommand{\A}{\mathcal{A}}
\newcommand{\B}{\mathcal{B}}
\newcommand{\cB}{\mathcal{B}}

\newcommand{\W}{\mathcal{W}}
\newcommand{\cP}{\mathscr{P}}

\newcommand{\pf}{\noindent{\bfseries Proof. }}

\newcommand{\bi}{\bs{\rm i}}
\newcommand{\bj}{\bs{\rm j}}

\newcommand{\M}{{\mathcal{M}}}

\newcommand{\gl}{\mathfrak{gl}}
\newcommand{\Z}{\mathbb{Z}}

\newcommand{\te}{\widetilde{e}}
\newcommand{\tf}{\widetilde{f}}

\numberwithin{equation}{section}

\begin{document}
\title[]
{A  plactic algebra of extremal weight crystals and the Cauchy
identity for Schur operators}
\author{JAE-HOON KWON}
\address{Department of Mathematics \\ University of Seoul   \\ Seoul 130-743, Korea }
\email{jhkwon@uos.ac.kr }

\begin{abstract}
We  give a new bijective interpretation of the Cauchy identity for
Schur operators which is a commutation relation between two formal
power series with operator coefficients. We introduce a plactic
algebra associated with the Kashiwara's extremal weight crystals
over the Kac-Moody algebra of type $A_{+\infty}$, and construct a
Knuth type correspondence preserving the plactic relations. This
bijection yields the Cauchy identity for Schur operators as a
homomorphic image of its associated identity for plactic characters
of extremal weight crystals, and also recovers the Sagan and
Stanley's correspondence for skew tableaux as its restriction.
\end{abstract} \maketitle

\section{Introduction}
Let $\Lambda=\Lambda_{{\bf x}}$ be the algebra of symmetric
functions  in formal commuting variables ${\bf
x}=\{\,x_1,x_2,\ldots\,\}$ over $\mathbb{Q}$. We denote by $\cP$ the
set of partitions and let $s_\lambda({\bf x})$ be the Schur function
in ${\bf x}$ corresponding to $\lambda\in\cP$. Let
\begin{equation*}
\begin{split}
\mathscr{P}({\bf x})=\sum_{\lambda\in\cP}s_\lambda s_{\lambda}({\bf x}),
 \ \ \mathscr{Q}({\bf x})=\sum_{\lambda\in\cP}s_\lambda^\perp s_{\lambda}({\bf x})\ \in {\rm End}_\mathbb{Q}(\Lambda)[[{\bf x}]],
\end{split}
\end{equation*}
where $s_{\lambda}$ and $s_{\lambda}^{\perp}$ are linear operators
on $\Lambda$ induced from the left multiplication by
$s_{\lambda}({\bf x})$ and its adjoint with respect to the Hall
inner product on $\Lambda$, respectively. One may regard $s_\lambda$
and $s_\lambda^\perp$ as operators on
$\mathbb{Q}\cP=\bigoplus_{\lambda\in\cP}\mathbb{Q}\lambda$, where
$\lambda$ is identified with $s_\lambda({\bf x})$. Moreover
$s_\lambda$ and $s_\lambda^\perp$ can be given as Schur functions in
certain locally non-commutative operators on $\mathbb{Q}\cP$ called
{\it Schur operators} by Fomin, while $\mathscr{P}({\bf x})$ and
$\mathscr{Q}({\bf y})$ can be written as Cauchy products in Schur
operators and ${\bf x}$ \cite{Fomin95,Fomin98}.

Let ${\bf y}=\{\,y_1,y_2,\ldots\,\}$ be another formal commuting
variables. It is well known that the following commutation relation
holds;
\begin{equation}\label{PQ-QP}
\mathscr{Q}({\bf y})\mathscr{P}({\bf x})=\mathscr{P}({\bf x})\mathscr{Q}({\bf y})
\frac{1}{\prod_{i,j}(1-x_iy_j)}
\end{equation}
called {\it generalized Cauchy identity} or {\it Cauchy identity for
Schur operators}. Considering both sides as operators with
coefficients in $\Lambda_{{\bf x}}\otimes \Lambda_{{\bf y}}$ and
then equating each entry of their matrix forms, we obtain a Cauchy
identity for skew Schur functions \cite{Mac95}
\begin{equation*}
\sum_{\lambda}s_{\lambda/\alpha}({\bf x})s_{\lambda/\beta}({\bf y})
=\sum_{\eta}s_{\beta/\eta}({\bf x})s_{\alpha/\eta}({\bf
y})\frac{1}{\prod_{i,j}(1-x_iy_j)},
\end{equation*}
where $\alpha, \beta$ are given partitions. A bijective
interpretation of the Cauchy identity for skew Schur functions was
given by Sagan and Stanley \cite{SS}, and it was extended to a
bijection in a more general framework by Fomin \cite{Fomin95}
including various analogues of Knuth correspondence.

Recently, a new representation theoretic interpretation of the
Cauchy identity for Schur operators was given by the author
\cite{K09} using the notion of Kashiwara's extremal weight crystals
\cite{Kas94'} over the quantized enveloping algebra associated with
the Kac-Moody algebra of type $A_{+\infty}$, say $\gl_{>0}$. It is
proved that a Schur operator can be realized as a functor of
tensoring by an extremal weight crystal element and (\ref{PQ-QP})
can be understood (but not in a bijective way) as a non-commutative
character identity corresponding to the decomposition of the crystal
graph of the Fock space with infinite positive level, which is an
infinite analogue of the level $n$ fermionic Fock space
decomposition due to Frenkel \cite{Fr}.

Motivated by a categorification of Schur operators in \cite{K09}, we
give a new combinatorial way to explain both the Cauchy identity for
Schur operators and skew Schur functions in terms of a single
bijection. More precisely, the main result in this paper is to
construct a Knuth type correspondence, which gives a bijective
interpretation of the identity (\ref{PQ-QP}) or its dual form, as
the usual Knuth correspondence does for the Cauchy product, and also
recovers Sagan and Stanley's correspondence as a restriction of this
bijection.

Our approach is to define a $t$-analogue of the plactic algebra
$\mathscr{U}(t)$ for $\gl_{>0}$ generated by $u_i$ and $u_{i^\vee}$
for $i\geq 1$ with $t$ an indeterminate, where the subalgebra
generated by $u_i$ (resp. $u_{i^\vee}$) is isomorphic to the usual
plactic algebra by Lascoux and Sch\"{u}tzenberger \cite{LS}. We show
that $\mathscr{U}(1)$ is isomorphic to the plactic algebra defined
by using the notion of crystal equivalence (cf.\cite{Lit}). Note
that each monomial in $\mathscr{U}(1)$ corresponds in general to an
element of an extremal weight crystal, which may not be either
highest or lowest weight crystal.

Now, let $\M_{\A,\B}$ be the set of $\A\times \B$ matrices
$A=(a_{ij})$ with entries in $\Z_{\geq 0}$ such that
$\sum_{i\in\A}\sum_{j\in\B}a_{ij}<\infty$ and $a_{ij}\leq 1$ for
$|i|\neq|j|$, where $\A$ and $\B$ are arbitrary $\Z_2$-graded sets
and $|\cdot|$ denotes the degree of an element in $\A$ or $\B$. We
assume that all the elements in $\mathbb{N}$ and $\mathbb{N}^\vee$
are of degree $0$. By using the usual Knuth map and non-commutative
Littlewood-Richardson rule  for extremal weight crystals for
$\gl_{>0}$ \cite{K09,K10}, we construct an explicit bijection
(Theorem \ref{main result});
$$\M_{\A,\mathbb{N}}\times \M_{\B,\mathbb{N}^\vee}
\longrightarrow \M_{\B,\mathbb{N}^\vee}\times \M_{\A,\mathbb{N}}
\times \M_{\A,\B},$$ which preserves the weights with respect to
$\A$ and $\B$, and the plactic relations of $\mathscr{U}(t)$ for the
column words with entries in $\mathbb{N}\cup\mathbb{N}^\vee$ on both
sides. As a corollary, we obtain a character identity in locally
non-commuting variables ${\bf u}=\{\,u_i,u_{i^\vee}\,|\,i\geq 1\,\}$
and commuting variables ${\bf x}_\A=\{\,x_a\,|\,a\in\A\,\}$, ${\bf
x}_\B=\{\,x_b\,|\,b\in\B\,\}$ (Corollary \ref{Cauchy identity}). In
particular, when $\A=\B=\mathbb{N}$, this identity recovers
(\ref{PQ-QP}) under a homomorphism sending $u_{i}$ and $u_{i^\vee}$
to Schur operators on $\mathbb{Q}\cP$ and specializing $t=1$.
Moreover, the Knuth correspondence for skew tableaux by Sagan and
Stanley can be recovered by restricting the above bijection to the
pairs of matrices on the left-hand side whose column words are
Littlewood-Richardson words of shape $(\alpha,\beta)$ with
$\alpha,\beta\in\cP$ (see Section \ref{skewRSK} for a definition).

The paper is organized as follow. In Section 2, we briefly recall
necessary background for semistandard tableaux and Knuth
correspondence. In Section 3, we recall the notion of rational
semistandard tableaux for $\gl_{>0}$ and their insertion algorithm.
In Section 4, we introduce a plactic algebra for $\gl_{>0}$
associated with rational semistandard tableaux. Finally, in Section
5, we construct a Knuth type correspondence and its associated
non-commutative character identity.

\section{Preliminary}
\subsection{Semistandard tableaux}
Throughout this paper, we assume that $\A$ (or $\B$) is a linearly
ordered $\mathbb{Z}_2$-graded set, that is, $\A=\A_0\sqcup\A_1$,
which is at most countable. We usually denote by $<$ a linear
ordering on a given linearly ordered $\mathbb{Z}_2$-graded set.  For
$a\in \A_\epsilon$ ($\epsilon\in\Z_2$),  we put $|a|=\epsilon$. By
convention, we let $\mathbb{N}=\{\,1<2<\cdots\,\}$, and
$[n]=\{\,1<\cdots<n\,\}$ for $n\geq 1$, where all the elements are
of degree $0$.

Let $\cP$ denote the set of partitions. We identify a partition
$\lambda=(\lambda_i)_{i\geq 1}$ with a Young diagram or a subset
$\{\,(i,j)\,|\,1\leq j\leq \lambda_i\,\}$ of
$\mathbb{N}\times\mathbb{N}$ following \cite{Mac95}. Let
$\ell(\lambda)=\left|\,\{\,i\,|\,\lambda_i\neq 0\,\}\,\right|$.
We denote by
$\lambda'=(\lambda'_i)_{i\geq 1}$ the conjugate partition of
$\lambda$ whose Young diagram is
$\{\,(i,j)\,|\,(j,i)\in\lambda\,\}$. For $\mu\in\cP$ with
$\lambda\supset\mu$, $\lambda/\mu$ denotes the skew Young diagram.

For a skew Young diagram $\lambda/\mu$, a tableau $T$ obtained by
filling $\lambda/\mu$ with entries in $\A$ is called
$\A$-semistandard if (1) the entries in each row (resp. column) are
weakly increasing from left to right (resp. from top to bottom), (2)
the entries in $\A_0$ (resp. $\A_1$) are strictly increasing in each
column (resp. row). We say that $\lambda/\mu$ is the shape of
$T$, and write ${\rm sh}(T)=\lambda/\mu$.

We denote by $SST_{\A}(\lambda/\mu)$ the set of all
$\A$-semistandard tableaux of shape $\lambda/\mu$. We set
$\cP_{\A}=\{\,\lambda\in\cP\,|\,SST_{\A}(\lambda)\neq
\emptyset\,\}$. For example, $\cP_{\A}=\cP$ when $\A$ is an infinite
set, and $\cP_{[n]}=\{\,\lambda\,|\,\ell(\lambda)\leq n\,\}$.

Let $\W_{\A}$ be the set of finite words with letters in $\A$. For
$T\in SST_{\A}(\lambda/\mu)$, we denote by $T(i,j)$ the entry of $T$
at $(i,j)\in \lambda/\mu$. We denote by $w_{\rm col}(T)\in \W_\A$
the word obtained by reading the entries of $T$ column by column
from right to left, and from top to bottom in each column. Also, we
denote by $w_{\rm row}(T)\in\W_\A$ the word obtained by reading the
entries of $T$ row by row from top to bottom, and from right to left
in each row.

Let $P_{\A}=\bigoplus_{a\in \A}\mathbb{Z}\epsilon_a$ be the free
abelian group with the basis $\{\,\epsilon_a\,|\,a\in \A\,\}$ and
let ${\bf x}_\A=\{\,x_a\,|\,a\in\A\,\}$ be a set of commuting formal
variables. For $\lambda=\sum_{a \in\A}\lambda_a\epsilon_a\in P_\A$,
let ${\bf x}_\A^\lambda=\prod_{a\in\A}x_a^{\lambda_a}$. For
$w=w_1\ldots w_r\in\W_\A$, we define ${\rm wt}_{\A}(w)=\sum_{1\leq
i\leq r}\epsilon_{w_i}\in P_{\A}$. For a skew Young diagram
$\lambda/\mu$ and $T\in SST_\A(\lambda/\mu)$, we define ${\rm
wt}_{\A}(T)=\sum_{(i,j)\in \lambda/\mu}\epsilon_{T(i,j)}$. Let
$s_{\lambda/\mu}({\bf x}_\A)=\sum_{T\in SST_{\A}(\lambda/\mu)}{\bf
x}_{\A}^{{\rm wt}_{\A}(T)}$, which is the character of
$SST_{\A}(\lambda/\mu)$. Note that $s_{\lambda/\mu}({\bf x}_\A)$ is
a usual (skew) Schur function when $\A=\mathbb{N}$.

We will also use the following operations on tableaux;

\begin{itemize}
\item[(1)] {\it dual} : Let $\A^{\vee}=\{\,a^\vee\,|\,a\in\A\,\}$ be the linearly ordered
$\mathbb{Z}_2$-graded set with $|a^\vee|=|a|$ and $a^\vee < b^\vee$
for $a>b$. For $T\in SST_{\A}(\lambda/\mu)$, we define $T^{\vee}$ to
be the tableau obtained by applying $180^{\circ}$-rotation to $T$
and replacing each entry $a$ in $T$ with $a^\vee$. Then $T^{\vee}\in
SST_{\A^{\vee}}((\lambda/\mu)^{\vee})$, where $(\lambda/\mu)^{\vee}$
denotes the shape of $T^{\vee}$. We use the convention that
$(a^\vee)^\vee=a$ for $a\in\A$ and hence $(T^\vee)^\vee=T$.

\item[(2)] {\it gluing} : Let $\A\ast\B$ be the $\mathbb{Z}_2$-graded set $\A\sqcup
\B$ with the extended linearly ordering given by $a<b$ for $a\in\A$
and $b\in\B$. For $S\in SST_{\A}(\mu)$ and $T\in
SST_{\B}(\lambda/\mu)$, we define $S\ast T \in
SST_{\A\ast\B}(\lambda)$ by $S\ast T(i,j)=S(i,j)$ for $(i,j)\in\mu$
and $T(i,j)$ for $(i,j)\in\lambda/\mu$.
\end{itemize}

\subsection{Littlewood-Richardson rule}\label{LRrule}
For $a\in\A$ and $T\in SST_{\A}(\lambda)$ with $\lambda\in\cP_{\A}$,
$a\rightarrow T$ (resp. $T\leftarrow a$) denotes the tableau
obtained by the Schensted column (resp. row) insertion (see for
example, \cite[Appendix A.2]{Fu} and \cite{BR} for its
super-analogue). For $w=w_1\ldots w_r\in\W_\A$, we let
$(w\rightarrow T)=( w_r\rightarrow(\cdots(w_1\rightarrow T)\cdots))$
and $(T \leftarrow w)=((\cdots(T\leftarrow w_1)\cdots )\leftarrow
w_r)$. For $S\in SST_{\A}(\mu)$ and $T\in SST_{\A}(\nu)$ with
$\mu,\nu\in \cP_{\A}$, we define $T\rightarrow S$ (resp.
$S\leftarrow T$) by $(w(T)_{\rm col}\rightarrow S)$ (resp.
$S\leftarrow \left(w(T)_{\rm row}\right)^{\rm rev}$), where $w^{\rm
rev}$ denotes the reverse word of $w$.

For $\lambda,\mu,\nu\in\cP$ with $|\lambda|=|\mu|+|\nu|$, let
${\rm\bf LR}^\lambda_{\mu \nu}$ be the set of tableaux $U$ in
$SST_{\mathbb{N}}(\lambda/\mu)$ such that
\begin{itemize}
\item[(1)] ${\rm wt}_\mathbb{N}(U)=\sum_{i\geq 1}\nu_i\epsilon_i$,

\item[(2)] for $1\leq k\leq |\nu|$, the number of occurrences of each $i\geq 1$ in
$w_1\ldots w_k$ is no less than that of $i+1$ in $w_1\ldots w_k$,
where $w(U)_{\rm col}=w_1\ldots w_{|\nu|}$.
\end{itemize}
We call ${\rm\bf LR}^\lambda_{\mu \nu}$ the set of
Littlewood-Richardson tableaux of shape $\lambda/\mu$ with content
$\nu$ and put $c^\lambda_{\mu \nu}=\left|{\rm\bf LR}^\lambda_{\mu
\nu}\right|$ \cite{Mac95}. Let us introduce a variation of ${\rm\bf
LR}^\lambda_{\mu \nu}$,  which is necessary for our later arguments.
Define $\overline{{\rm\bf LR}}^\lambda_{\mu \nu}$ to be the set of
tableaux $U$ in $SST_{-\mathbb{N}}(\lambda/\mu)$ such that
\begin{itemize}
\item[(1)] ${\rm wt}_{-\mathbb{N}}(U)=\sum_{i\geq 1}\nu_i\epsilon_{-i}$,

\item[(2)] for $1\leq k\leq |\nu|$, the number of occurrences of each $-i\leq -1$ in
$w_k\ldots w_{|\nu|}$ is no less than that of $-(i+1)$ in $w_k\ldots
w_{|\nu|}$, where $w(U)_{\rm col}=w_1\ldots w_{|\nu|}$.
\end{itemize}

Note that for $U\in SST_{\mathbb{N}}(\lambda/\mu)$, $U\in {\rm\bf
LR}^\lambda_{\mu \nu}$ if and only if $U$ is Knuth equivalent to
$H_{\nu}\in SST_{\mathbb{N}}(\nu)$, where $H_\nu(i,j)=i$ for
$(i,j)\in \nu$ (cf.\cite{Fu}). Similarly, we have for $U\in
SST_{-\mathbb{N}}(\lambda/\mu)$, $U\in \overline{{\rm\bf
LR}}^\lambda_{\mu \nu}$ if and only if $U$ is Knuth equivalent to
$L_{\nu} \in SST_{-\mathbb{N}}(\nu)$, where $L_\nu(i,j)=-\nu'_j+i-1$
for $(i,j)\in\nu$.

There is also a one-to-one correspondence from the set of $V\in
SST_{\mathbb{N}}(\nu)$ such that $(V \rightarrow H_\mu) = H_\lambda$
to ${\rm\bf LR}^\lambda_{\mu \nu}$. Indeed, $V$ corresponds to
$\imath(V)\in {\rm\bf LR}^\lambda_{\mu \nu}$ where the number of
$k$'s in the $i$-th row of $V$ is equal to the number of $i$'s in
the $k$-th row of $\imath(V)$ for $i, k\geq 1$.

\begin{ex}{\rm\mbox{}
\begin{center}
$SST_{\mathbb{N}}((3,3,2)) \ni$
\begin{tabular}{ccc}
    1 & 1 & 2  \\
    2 & 2 & 3  \\
    3 & 4
\end{tabular} \ \
\ \ $\stackrel{\imath}{\longrightarrow}$ \ \ \ \
\begin{tabular}{ccccc}
    {$\bullet$} & {$\bullet$} & {$\bullet$} &  1 & 1 \\
    {$\bullet$} & 1 & 2 & 2 & \\
    2 & 3 & & & \\
    3 &  & & &
\end{tabular}
$\in {\bf LR}^{(5,4,2,1)}_{(3,1)\ (3,3,2)}$
\end{center}\vskip 2mm}
\end{ex}

For $S\in SST_{\A}(\mu)$ and $T\in SST_{\A}(\nu)$, suppose that
$\lambda={\rm sh}(T\rightarrow S)$ and $w_{\rm col}(T)=w_1\ldots
w_r$. Define $(T\rightarrow S)_R$ to be the tableau of shape
$\lambda/\mu$ such that ${\rm sh}(w_1\ldots w_k\rightarrow S)/{\rm
sh}(w_1\ldots w_{k-1}\rightarrow S)$ is filled with $i$ if ${w_k}$
is in the $i$th row of $T$ for $1\leq k\leq r$. Then the map
$(S,T)\mapsto \left((T\rightarrow S),(T\rightarrow S)_R\right)$
gives a bijection \cite{Thomas}
\begin{equation}\label{LR}
SST_{\A}(\mu)\times SST_{\A}(\nu)\longrightarrow
\bigsqcup_{\lambda\in\cP_{\A}}SST_{\A}(\lambda)\times {\bf
LR}^{\lambda}_{\mu\,\nu},
\end{equation}
which also implies $s_\mu({\bf x}_{\A})s_\nu({\bf
x}_{\A})=\sum_\lambda c^\lambda_{\mu \nu}s_\lambda({\bf x}_{\A})$.

\subsection{Skew Littlewood-Richardson rule}\label{skewLRrule}
Let $\lambda/\mu$ be a skew Young diagram. Let $U$ be a tableau of
shape $\lambda/\mu$ with entries in $\mathcal{A}\sqcup\mathcal{B}$,
satisfying the following conditions;
\begin{itemize}
\item[(S1)] $U(i,j)\leq U(i',j')$ whenever $U(i,j), U(i',j')\in\mathcal{X}$
for $(i,j), (i',j')\in\lambda/\mu$ with $i\leq i'$ and $j\leq j'$,

\item[(S2)] in each column of $U$,
entries in $\mathcal{X}_0$  increase strictly from top to bottom,

\item[(S3)] in each row of $U$,
entries in $\mathcal{X}_1$  increase strictly from left to right,
\end{itemize}
where  $\mathcal{X}=\A$ or $\mathcal{B}$. Suppose that $a\in
\mathcal{A}$ and $b\in \mathcal{B}$ are two adjacent entries in $U$
such that $b$ is placed above or to the left of $a$. Interchanging
$a$ and $b$ is called a switching if the resulting tableau still
satisfies the conditions (S1), (S2) and (S3).

For $T\in SST_{\mathcal{A}}(\lambda/\mu)$, let $U$ be a tableau
obtained from $H_\mu\ast T$ by applying switching procedures as far
as possible (in this case, $\mathcal{B}=\mathbb{N}$). Then
$U=\jmath(T)\ast \jmath(T)_R$ for some $\jmath(T)\in SST_\A(\nu)$
and $\jmath(T)_R\in SST_{\mathbb{N}}(\lambda/\nu)$ with $\nu\in\cP$.
Then by \cite[Theorem 3.1]{BSS}, the map sending $T$ to
$(\jmath(T),\jmath(T)_R)$ gives a bijection
\begin{equation}\label{skewLR}
SST_{\A}(\lambda/\mu) \longrightarrow \bigsqcup_{\nu\in\cP_\A}
SST_{\A}(\nu)\times {\bf LR}^{\lambda}_{\nu\,\mu}.
\end{equation}
In particular, the map $Q \mapsto \jmath(Q)_R$ restricts to a
bijection from ${\bf LR}^{\lambda}_{\mu\,\nu}$ to ${\bf
LR}^{\lambda}_{\nu\,\mu}$, and from $\overline{\bf
LR}^{\lambda}_{\mu\,\nu}$ to ${\bf LR}^{\lambda}_{\nu\,\mu}$ when
$\A=\pm\mathbb{N}$, respectively.

\subsection{RSK correspondence}

Let $\M_{\A,\B}$ be the set of $\A\times \B$ matrices $A=(a_{ij})$
with entries in $\Z_{\geq 0}$ such that
$\sum_{i\in\A}\sum_{j\in\B}a_{ij}<\infty$ and $a_{ij}\leq 1$ for
$|i|\neq|j|$. Let $\Omega_{\A,\B}$ be the set of biwords
$(\bi,\bj)\in \W_{\A}\times \W_{\B}$ such that
\begin{itemize}
\item[(1)] $\bi=i_1\cdots i_r$ and $\bj=j_1\cdots j_r$  for some $r\geq
0$,

\item[(2)] $(i_1,j_1)\leq \cdots \leq (i_r,j_r)$,

\item[(3)] $(i_s,j_s)<
(i_{s+1},j_{s + 1})$ if $|i_s|\neq|j_s|$ for $1\leq s<r$,
\end{itemize}
where for $(i,j)$ and $(k,l)\in \A\times \B$,
\begin{equation*}\label{partial order}
(i,j)< (k,l) \ \ \ \ \Longleftrightarrow \ \ \ \
\begin{cases}
(i<k) & \text{or}, \\
(i=k,\ |i|=0,\ \text{and} \ j>l) & \text{or}, \\
(i=k, \ |i|=1,\ \text{and} \ j<l) &.
\end{cases}
\end{equation*}
There is a bijection from $\Omega_{\A,\B}$ to $\M_{\A,\B}$, where
$(\bi,\bj)$ is mapped to $A(\bi,\bj)=(a_{ij})$ with
$a_{ij}=|\{\,k\,|\,(i_k,j_k)=(i,j) \,\}|$. Note that the pair of
empty words $(\emptyset,\emptyset)$ corresponds to zero matrix.

For $A=A(\bi,\bj)\in \M_{\A,\B}$, we let $P(A)=(\bj \rightarrow
\emptyset)$  and let $Q(A)$ be the tableau of the same shape as
$P(A)$ such that ${\rm sh}(j_1\ldots j_k \rightarrow \emptyset)/{\rm
sh}(j_1\ldots j_{k-1} \rightarrow \emptyset)$ is filled with $i_k$
for $k\geq 1$. Then the map sending $A$ to $(P(A),Q(A))$ gives a
bijection
\begin{equation}\label{Knuth-1}
\M_{\A,\B} \longrightarrow \bigsqcup_{\lambda\in
\cP_\A\cap\cP_{\B}}SST_\B(\lambda)\times SST_{\A}(\lambda),
\end{equation}
which is the (super analogue of) RSK correspondence \cite{Kn}. If we
define ${\rm wt}_\A(A)={\rm wt}_\A(\bj)$ and ${\rm wt}_\B(A)={\rm
wt}_\A(\bi)$, then the bijection preserves ${\rm wt}_\A$ and ${\rm
wt}_\B$. In terms of characters, we obtain the following Cauchy
identity
\begin{equation*}
\frac{\prod_{|a|\neq |b|}(1+x_ax_b)}{\prod_{|a|=
|b|}(1-x_ax_b)}=\sum_{\lambda\in \cP_\A\cap\cP_{\B}}s_{\lambda}({\bf
x}_\B)s_\lambda({\bf x}_\A),
\end{equation*}
where $a\in\A$ and $b\in\B$.

Similarly, for $A=(a_{ij})\in \M_{\A,\B}$, let $A'=(a'_{ij^\vee})$
be the unique matrix in $\M_{\A,\B^\vee}$ such that
$a_{ij}=a'_{i,j^\vee}$ for $(i,j)\in\A\times\B$. Then the map
sending $A$ to $(P(A')^\vee,Q(A'))$ gives a bijection
\begin{equation}\label{Knuth-2}
\M_{\A,\B} \longrightarrow \bigsqcup_{\lambda\in
\cP_\A\cap\cP_{\B}}SST_\B(\lambda^\vee)\times SST_{\A}(\lambda).
\end{equation}

Finally, for $\mu\in\cP_\B$, we have {\allowdisplaybreaks
\begin{align*}
 SST_\B(\mu)\times
&\M_{\A,\B}    \\
&\stackrel{1-1}{\longleftrightarrow} \bigsqcup_{\mu,\in\cP_\B
\nu\in\cP_\A}SST_\B(\mu)\times
SST_{\B}(\nu)\times SST_{\A}(\nu)\ \ \ \ \ \ \ \ \ \text{by (\ref{Knuth-1})}\\
&\stackrel{1-1}{\longleftrightarrow} \bigsqcup_{\substack{\lambda\in
\cP_\A\cap\cP_{\B} \\ \mu\subset \lambda}}SST_\B(\lambda)\times
\left(\bigsqcup_{\nu\in\cP_\A} SST_{\A}(\nu)\times {\bf
LR}^{\lambda}_{\nu\,\mu}\right)  \ \ \text{by (\ref{LR})} \\
&\stackrel{1-1}{\longleftrightarrow}\ \
\bigsqcup_{\substack{\lambda\in \cP_\A\cap\cP_{\B} \\ \mu\subset
\lambda}}SST_\B(\lambda)\times SST_{\A}(\lambda/\mu)\ \ \ \ \ \ \ \ \ \ \ \ \ \ \ \ \ \ \ \ \ \text{by (\ref{skewLR})}\\
\end{align*}}
Hence we obtain a bijection
\begin{equation}\label{Knuth-3}
SST_\B(\mu)\times \M_{\A,\B} \longrightarrow
\bigsqcup_{\substack{\lambda\in \cP_\A\cap\cP_{\B} \\ \mu\subset
\lambda}}SST_\B(\lambda)\times SST_{\A}(\lambda/\mu).
\end{equation}

\section{Rational semistandard tableaux}
\subsection{Rational semistandard
tableaux for $\gl_{>0}$}\label{insertion} For convenience, we let
for a skew Young diagram $\lambda/\mu$,
$$\B_{\lambda/\mu}=SST_{\mathbb{N}}(\lambda/\mu),
\ \ \ \ \ \B_{\lambda/\mu}^\vee= SST_{\mathbb{N}^\vee}((\lambda/\mu)^\vee).$$

\begin{df}\label{Bmunu}{\rm
For $\mu,\nu\in \cP$, we define $\cB_{\mu,\nu}$ to be the set of
bitableaux $(S,T)$ such that
\begin{itemize}
\item[(1)] $(S,T)\in \B_{\mu}\times \B_{\nu}^\vee$,

\item[(2)]
$\Bigl|\,\{\,i\,|\,S(i,1)\leq k\,\}\,
\Bigr|+\left|\,\{\,i\,|\,T(i,1)\geq k^\vee\,\}\, \right| \leq k$ for
$k\geq 1$.
\end{itemize}
}
\end{df}
For convenience, we identify $\B_{\mu,\emptyset}$ and
$\B_{\emptyset,\nu}$ with $\B_\mu$ and $\B^\vee_\nu$, respectively.

\begin{ex}{\rm
$$
\left(\begin{array}{ccc}
1 & 1   & 3      \\
2 & 3   &        \\
4 &     &
\end{array}, \ \
\begin{array}{ccc}
     &  7^\vee  \\
    5^\vee &  5^\vee  \\
     4^\vee  & 3^\vee
\end{array}
\right)\in \cB_{(3,2,1),(2,2,1)}$$}
\end{ex}

\begin{rem}{\rm
 Let $\gl_{>0}$ be the general linear Lie
algebra spanned by $\mathbb{N}\times\mathbb{N}$ complex matrices of
finite support. Then $\cB_{\mu,\nu}$ parameterizes a basis of an
extremal weight module over the quantum group $U_q(\gl_{>0})$
\cite{K09}. Recall that $\cB_{\mu,\nu}\cap
\Bigl(SST_{[n]}(\mu)\times SST_{[n]^\vee}(\nu^\vee)\Bigr)$ ($n\geq
2$) parameterizes a basis of a finite dimensional complex
irreducible representation of $\gl_n$  \cite{St}.}
\end{rem}

Let us review an insertion algorithm for $\B_{\mu,\nu}$ \cite{K09},
which is an infinite analogue of those for rational semistandard
tableaux for $\gl_n$ \cite{St,Str}. For $a\in\mathbb{N}$ and
$(S,T)\in\cB_{\mu,\nu}$, we define $a\rightarrow (S,T)$ in the
following way;

Suppose first that $S$ is empty and $T$ is a single column tableau.
Let $(T',a')$ be the pair obtained as follows.

\begin{itemize}
\item[(1)] If $T$ contains $a^\vee, (a+1)^\vee,\ldots, (b-1)^\vee$
but not $b^\vee$, then $T'$ is the tableau obtained from $T$ by
replacing $a^\vee, (a+1)^\vee,\ldots, (b-1)^\vee$ with
$(a+1)^\vee,(a+2)^\vee,\ldots, b^\vee$, and put $a'=b$.

\item[(2)] If $T$ does not contain $a^\vee$, then leave $T$
unchanged and put $a'=a$.
\end{itemize}
Now, we suppose that $S$ and $T$ are arbitrary.
\begin{itemize}
\item[(1)] Apply the above process to the leftmost column of $T$ with $a$.

\item[(2)] Repeat (1) with $a'$  and the next column to the right.

\item[(3)] Continue this process to the right-most column of $T$ to get a tableau $T'$ and $a''$.

\item[(4)] Define
$$(a\rightarrow (S,T))=\left((a''\rightarrow S)\,,\,T'\right).$$
\end{itemize}
Then $a\rightarrow (S,T)\in\cB_{\sigma,\nu}$ for some $\sigma
\in\cP$ with $|\sigma/\mu|=1$. For $w=w_1\ldots w_r\in \W_
\mathbb{N}$, we let $(w\rightarrow (S,T))=(
w_r\rightarrow(\cdots(w_1\rightarrow (S,T))\cdots))$.\vskip 3mm

Next, we define $(S,T)\leftarrow a^\vee$ to be the pair $(S',T')$
obtained in the following way;
\begin{itemize}
\item[(1)] If the pair $(S,(T^\vee \leftarrow a)^\vee)$ satisfies the condition (2) in Definition \ref{Bmunu},
then put $S'=S$ and $T'=(T^\vee \leftarrow a)^\vee$.

\item[(2)] Otherwise, choose the smallest $k$
such that $a_k$ is bumped out of the $k$-th row in the row insertion
of $a$ into $T^\vee$ and the insertion of $a_k$ into the $(k+1)$-st
row violates the condition (2) in Definition \ref{Bmunu}.

\item[(2-a)] Stop the row insertion of $a$ into
$T^\vee$ when $a_k$ is bumped out and let $T'$ be the resulting
tableau after taking $\vee$.

\item[(2-b)] Remove $a_k$ in the left-most column of $S$, which
necessarily exists, and then play the jeu de taquin  (see for
example \cite[Section 1.2]{Fu}) to obtain a tableau $S'$.
\end{itemize}
In this case, $(S,T)\leftarrow a^\vee\in \cB_{\sigma,\tau}$, where
either (1) $|\mu/\sigma|=1$ and $\tau=\nu$, or (2) $\sigma=\mu$ and
$|\tau/\nu|=1$. For $w=w_1\ldots w_r\in \W_{\mathbb{N}^\vee}$, we
let $((S,T)\leftarrow w)=((\cdots((S,T)\leftarrow w_1)\cdots
)\leftarrow w_r)$.\vskip 3mm

\subsection{Non-commutative Littlewood-Richardson rule}
Let us recall the Littlewood-Richardson rule  for $\B_{\mu,\nu}$
(see \cite[Proposition 4.9]{K09} for more details).

Suppose that $\mu,\nu \in\cP$ are given. For $(S,T)\in
\B_\nu^\vee\times \B_\mu$,  consider $(w_{\rm col}(T)\rightarrow
(\emptyset,S))$. Suppose that $w_{\rm col}(T)=w_1\ldots w_r$ and
$(w_1\ldots w_k\rightarrow (\emptyset,S))\in\B_{\mu^{(k)},\nu}$ for
$1\leq k\leq r$ with $\mu^{(r)}=\mu$. Let $(i_k,j_k)\in \mu$
correspond to $w_k$ in $T$ ($1\leq k\leq r$). Then $\mu^{(k)}$ is
obtained from $\mu^{(k-1)}$ by adding a box in the $i_k$-th row.
Hence,  the map sending $(S,T)$ to $(w_{\rm col}(T)\rightarrow
(\emptyset,S))$ gives a bijection \cite[Corollary 4.11]{K09}
\begin{equation}\label{bijection-1}
\B_\nu^\vee \times \B_{\mu} \longrightarrow \B_{\mu,\nu}.
\end{equation}
\vskip 3mm

Next, let $(S,T)\in \B_\mu\times \B_\nu^\vee$ be given.  Consider
$((S,\emptyset)\leftarrow w_{\rm col}(T))$. Suppose that $w_{\rm
col}(T)=w_1\ldots w_r$ and $((S,\emptyset)\leftarrow w_1\ldots
w_k)\in\B_{\mu^{(k)},\nu^{(k)}}$ for $1\leq k\leq r$. Let
$(i_k,j_k)\in \nu$ correspond to $w_k$ in $T$ ($1\leq k\leq r$).
Define $U$ to be the tableau of shape $\nu$ such that for $1\leq
k\leq r$
\begin{equation*}
\begin{split}
&U(i_k,j_k) \\ &=\begin{cases} i, & \text{if $\mu^{(k)}$ is obtained
from $\mu^{(k-1)}$  by removing a box in the $i$-th row},
\\ -j, & \text{if $\nu^{(k)}$ is
obtained from $\nu^{(k-1)}$  by adding a box in the $j$-th row}.
\end{cases}
\end{split}
\end{equation*}
Then $U=U_+\ast U_-$, where $U_{+}\in SST_{\mathbb{N}}(\lambda)$ and
$U_{-}\in SST_{-\mathbb{N}}(\nu/\lambda)$ for some $\lambda\subset
\nu$. Let $\sigma=\mu^{(r)}$ and $\tau=\nu^{(r)}$. We have
$\imath(U_{+})\in {\bf LR}^{\mu}_{\sigma \lambda}$ and $U_{-}\in
\overline{\bf LR}^{\nu}_{\lambda \tau}$ for some $\lambda$, hence
$\jmath(U_{-})_R\in {\bf LR}^{\nu}_{\tau\lambda }$ (see Section
\ref{LRrule}, \ref{skewLRrule}). Therefore, we have a bijection
\cite[Proposition 4.3]{K10}
\begin{equation}\label{bijection-2}
\cB_\mu\times \cB_\nu^\vee {\longrightarrow}
\bigsqcup_{\lambda,\sigma,\tau}\cB_{\sigma,\tau}\times{\bf
LR}^\mu_{\sigma \lambda}\times {\bf LR}^{\nu}_{\tau \lambda},
\end{equation}
where $(S, T)$ is mapped to
$\left(((S,\emptyset)\leftarrow w_{\rm col}(T)),\imath(U_{+}), \jmath(U_{-})_R \right)$.

Now, we have {\allowdisplaybreaks
\begin{align*}
 \bigsqcup_{\lambda,\sigma,\tau }&\cB_{\sigma,\tau}\times{\bf
LR}^\mu_{\sigma \lambda}\times {\bf LR}^{\nu}_{\tau \lambda}\\
&\stackrel{1-1}{\longleftrightarrow} \bigsqcup_{\lambda,\sigma,\tau
}\B_{\tau}^\vee\times \B_\sigma\times{\bf
LR}^\mu_{\sigma \lambda}\times {\bf LR}^{\nu}_{\tau \lambda} \ \ \ \ \ \text{by (\ref{bijection-1})} \\
&\stackrel{1-1}{\longleftrightarrow} \bigsqcup_{\lambda,\sigma,\tau
}\B_{\tau}^\vee\times {\bf LR}^{\nu}_{\tau \lambda} \times
\B_\sigma\times {\bf
LR}^\mu_{\sigma \lambda} \\
&\stackrel{1-1}{\longleftrightarrow}\ \bigsqcup_{\lambda\subset
\mu,\nu}\cB_{\nu/\lambda}^\vee\times \cB_{\mu/\lambda}  \ \ \ \ \  \
\ \ \ \ \ \ \ \ \ \ \ \ \ \ \text{by (\ref{skewLR})}.
\end{align*}}
Hence, we obtain the following bijection \cite[Proposition 5.1]{K10}
\begin{equation}\label{extLR}
\cB_\mu\times \cB_\nu^\vee \longrightarrow \bigsqcup_{\lambda\subset
\mu,\nu}\cB_{\nu/\lambda}^\vee\times \cB_{\mu/\lambda}.
\end{equation}

\section{Plactic algebra}
\subsection{A plactic algebra for $\gl_{>0}$}
Let $t$ be an indeterminate. Define $\mathscr{U}(t)$ to be an
associative $\mathbb{Q}[t,t^{-1}]$-algebra with unity generated by
$u_i$ and $u_{i^\vee}$ ($i\in\mathbb{N}$) subject to the following
relations;
\begin{equation}\label{plactic}
\begin{split}
&u_iu_ju_k=u_iu_ku_j, \ \ u_{k^\vee}u_{j^\vee}u_{i^\vee}=u_{j^\vee}u_{k^\vee}u_{i^\vee} \ \ (j\leq i<k),\\
&u_iu_ju_k=u_ju_iu_k, \ \ u_{k^\vee}u_{j^\vee}u_{i^\vee}=u_{k^\vee}u_{i^\vee}u_{j^\vee} \ \ (j< i\leq k),\\
& u_{i+1}u_{(i+1)^\vee}=u_{i^\vee}u_i \ \ \ (i\geq 1), \ \ \
u_1u_{1^\vee}=t,\\
& u_iu_{j^\vee}=u_{j^\vee}u_i \ \ \ (i\neq j).
\end{split}
\end{equation}
Let $\mathscr{U}(t)_+$ (resp. $\mathscr{U}(t)_-$) be the subalgebra
of $\mathscr{U}(t)$ generated by $u_{i}$ (resp. $u_{i^\vee}$) for
$i\geq 1$. Then $\mathscr{U}(t)_{\pm}$ is isomorphic to the usual
plactic algebra for $\gl_{>0}$ over $\mathbb{Q}[t,t^{-1}]$
\cite{LS}, where the first two relations in (\ref{plactic}) are
Knuth relations.

Let $\W$ be the set of finite words with letters in
${\mathbb{N}\cup\mathbb{N}^\vee}$. For $w=w_1\cdots w_r\in \W$, put
$u_{w}=u_{w_1}\cdots u_{w_r}\in \mathscr{U}(t)$. It is well-known
that if $w \in\W_\mathbb{N}$ (resp. $\W_{\mathbb{N}^\vee}$), then
there exists a unique $T\in \B_\mu$ (resp. $\B_\mu^\vee$) for some
$\mu\in\cP$ such that $u_w=u_{w_{\rm col}(T)}=u_{w_{\rm row}(T)}$.

For a skew Young diagram $\lambda/\mu$ and $T\in \B_{\lambda/\mu}$
or $\B_{\lambda/\mu}^\vee$, we let $u_T=u_{w_{\rm col}(T)}$, and for
$\mu,\nu\in \cP$ and $(S,T)\in\B_{\mu,\nu}$, we let
$u_{(S,T)}=u_{S}u_{T}$.

\begin{lem}\label{reduction-0}
For $p,q\geq 1$, let $S\in \B_{(1^p)}$ and $T\in\B_{(1^q)}^\vee$ be
given and let $(S',T')=(w_{\rm col}(S)\rightarrow T)\in
\B_{(1^p),(1^q)}$. Then $u_{T}u_S=u_{S'}u_{T'}$.
\end{lem}
\pf It is straightforward to check from (\ref{bijection-1}) and
(\ref{plactic}). \qed

\begin{lem}\label{reduction}
For $p,q\geq 1$, let $S\in \B_{(1^p)}$ and $T\in\B_{(1^q)}^\vee$ be
given with $w_{\rm col}(S)=w^+_1\ldots w^+_p$ and $w_{\rm
col}(T)=w^-_q\ldots w^-_1$. Suppose that there exists $k\geq 1$ such
that $\Bigl|\,\{\,i\,|\,w^+_i\leq k\,\}\,
\Bigr|+\left|\,\{\,j\,|\,w^-_j\geq k^\vee\,\}\, \right| > k$. If
$w^+_i=k$ and $w^-_j=k^\vee$ for some $i$ and $j$, and $(S,T')\in
\B_{(1^p),(1^{q-1})}$, where $T'$ is obtained from $T$ by removing
$k^\vee$, then
$$u_{S}u_{T}=t\,u_{w^+_1\ldots\widehat{{w^+_i}}\ldots
w^+_p}u_{w^-_q \ldots\widehat{w^-_j}\ldots w^-_1}.$$
\end{lem}
\pf We use induction on $p+q$. If $p+q=2$, then $k=1$ and
$u_{w^+_1}u_{w^-_1}=u_1u_{1^\vee}=t$.

Suppose that $p+q\geq 2$ with $i<p$ and $j<q$. Note that
\begin{equation*}
\begin{split}
w^+_{i+1}\ldots w^+_{p}&=(k+a_1)\ldots (k+a_{p-i}), \\
w^-_q\ldots w^-_{j+1}&=(k+b_{q-j})^\vee\ldots (k+b_{1})^\vee,
\end{split}
\end{equation*}
for some $a_1<\cdots <a_{p-i}$ and $b_1<\cdots < b_{q-j}$. Also it
follows from our hypothesis that $i+j=k+1$, and hence we can choose
$(\overline{S},\overline{T})\in \B_{(1^{p-i}),(1^{q-j})}$ such that
$w_{\rm col}(\overline{S})=a_1 \ldots a_{p-i}$ and $w_{\rm
col}(\overline{T})=b_{q-j}^\vee\ldots b_{1}^\vee$. By Lemma
\ref{reduction-0}, there exists $(\overline{T}',\overline{S}')\in
\B_{(1^{q-j})}^\vee\times\B_{(1^{p-i})}$ with $w_{\rm
col}(\overline{S}')=a'_1\ldots a'_{p-i}$ and $w_{\rm
col}(\overline{T}')=(b'_{q-j})^\vee\ldots (b'_1)^\vee$ such that
$u_{\overline{T}'}u_{\overline{S}'}=u_{\overline{S}}u_{\overline{T}}$.
This implies that
\begin{equation*}
\begin{split}
&u_{w^+_{i+1}\ldots w^+_{p}}u_{w^-_q\ldots w^-_{j+1}} =
u_{(k+b'_{q-j})^\vee}\ldots u_{(k+b'_{1})^\vee}u_{(k+a'_1)}\ldots
u_{(k+a'_{p-i})}
\end{split}
\end{equation*}

Since $w^+_i=k< k+b'_{1}$ and $w^-_j=k^\vee > (k+a'_1)^\vee$, we
have
\begin{equation*}
u_{S}u_{T}=u_{(k+b'_{q-j})^\vee\ldots
(k+b'_{1})^\vee}\left(u_{w^+_1}\ldots{u_{w^+_i}}u_{w^-_j}\ldots
u_{w^-_1}\right)u_{(k+a'_1)\ldots (k+a'_{p-i})},
\end{equation*}
and by induction hypothesis,
\begin{equation*}
\begin{split}
u_{S}u_{T}&= t\,u_{(k+b'_{q-j})^\vee\ldots
(k+b'_{1})^\vee}\left(u_{w^+_1}\ldots\widehat{u_{w^+_i}}\widehat{u_{w^-_j}}\ldots
u_{w^-_1}\right)u_{(k+a'_1)\ldots (k+a'_{p-i})} \\
&=t\,u_{w^+_1}\ldots\widehat{u_{w^+_i}}\ldots
u_{w^+_p}u_{w^-_q}\ldots\widehat{u_{w^-_j}}\ldots u_{w^-_1}.
\end{split}
\end{equation*}
Now, we assume that $i=p$ and $j=q$, that is, $w^+_p=k$ and
$w^-_q=k^\vee$. Note that $p+q=k+1$. If $w^+_{p-1}\neq k-1$ and
$w^-_{q-1}\neq (k-1)^\vee$, then
$$\Bigl|\,\{\,i\,|\,w^+_i\leq k-2\,\}\,
\Bigr|+\left|\,\{\,j\,|\,w^-_j\geq (k-2)^\vee\,\}\, \right| =
p+q-2=k-1>k-2,$$ which contradicts the fact that $(S,T')\in
\B_{(1^p),(1^{q-1})}$. So we also assume that either $w^+_{p-1}=
k-1$ or $w^-_{q-1}= (k-1)^\vee$. \vskip 2mm

\noindent \textsc{Case 1}. Suppose that $w^+_{p-1}\neq k-1$ and
$w^-_{q-1}= (k-1)^\vee$. We have {\allowdisplaybreaks
\begin{align*}
u_{S}u_{T}&=u_{w^+_1}\ldots u_{w^+_{p-1}}u_k u_{k^\vee}
u_{w^-_{q-1}}\ldots u_{w^-_1}\\
&=u_{w^+_1}\ldots u_{w^+_{p-1}}u_{(k-1)^\vee}u_{k-1}
u_{w^-_{q-1}}\ldots u_{w^-_1}\\
&=u_{(k-1)^\vee} u_{w^+_1}\ldots u_{w^+_{p-1}}u_{k-1}u_{(k-1)^\vee}
u_{w^-_{q-2}}\ldots u_{w^-_1}\\
&=t\,u_{(k-1)^\vee} u_{w^+_1}\ldots u_{w^+_{p-1}}
u_{w^-_{q-2}}\ldots u_{w^-_1}\\
&=t\,u_{w^+_1}\ldots u_{w^+_{p-1}}
u_{(k-1)^\vee} u_{w^-_{q-2}}\ldots u_{w^-_1}\\
&=t\,u_{w^+_1}\ldots{u_{w^+_{p-1}}}{u_{w^-_{q-1}}}\ldots u_{w^-_1},
\end{align*}}
where we use induction hypothesis in the third line.

\noindent \textsc{Case 2}. Suppose that $w^+_{p-1}= k-1$ and
$w^-_{q-1}\neq (k-1)^\vee$. By almost the same argument as in
\textsc{Case 1}, we have
$$u_{S}u_{T}=t\,u_{w^+_1}\ldots{u_{w^+_{p-1}}}{u_{w^-_{q-1}}}\ldots
u_{w^-_1}.$$

\noindent \textsc{Case 3}. Suppose that $w^+_{p-1}= k-1$ and
$w^-_{q-1}= (k-1)^\vee$. We have {\allowdisplaybreaks
\begin{align*}
u_{S}u_{T}&=u_{w^+_1}\ldots u_{w^+_{p-1}}u_k u_{k^\vee}
u_{w^-_{q-1}}\ldots u_{w^-_1}\\
&=u_{w^+_1}\ldots u_{w^+_{p-1}}u_{(k-1)^\vee}u_{k-1}
u_{w^-_{q-1}}\ldots u_{w^-_1}\\
&=u_{(k-a)^\vee} u_{v_1}\ldots u_{v_{p-2}}u_{k-2}u_{k-1}
u_{w^-_{q-1}}\ldots u_{w^-_1},
\end{align*}}
for some $1 \leq a< k$ and $1\leq v_1<\ldots < v_{p-2}\leq k-2$. So
by induction hypothesis, {\allowdisplaybreaks
\begin{align*}
u_{S}u_{T}&=u_{(k-a)^\vee} u_{v_1}\ldots
u_{v_{p-2}}u_{k-2}u_{k-1}u_{(k-1)^\vee} u_{w^-_{q-2}}\ldots u_{w^-_1} \\
&=t\,u_{(k-a)^\vee} u_{v_1}\ldots
u_{v_{p-2}}u_{k-2}  u_{w^-_{q-2}}\ldots u_{w^-_1} \\
&=t\,u_{w^+_1}\ldots{u_{w^+_{p-2}}}u_{k-1}u_{(k-1)^\vee}
{u_{w^-_{q-1}}}\ldots
u_{w^-_1}\\
&=t\,u_{w^+_1}\ldots{u_{w^+_{p-1}}}{u_{w^-_{q-1}}}\ldots u_{w^-_1}.
\end{align*}}
This completes the induction. \qed

\begin{lem}\label{plactic lemma} Let $\mu,\nu\in\cP$ be given. For $a\in \mathbb{N}$ and $(S,T)\in\B_{\mu,\nu}$,
\begin{itemize}
\item[(1)] $u_{(S,T)}u_a=u_{(a\rightarrow(S,T))}$,

\item[(2)] $u_{(S,T)}u_{a^\vee}=t^\epsilon\,u_{((S,T)\leftarrow
a^\vee)}$, where $\epsilon=0,1$.
\end{itemize}
\end{lem}
\pf We keep the notations in Section \ref{insertion}. Consider
$(a\rightarrow(S,T))=(S',T')$. Let $(T',a')$ be the pair obtained by
the first step in the definition of $a\rightarrow (S,T)$. It is
straightforward to check that $u_{T}u_a = u_{a'}u_{T'}$. Since
$S'=(a'\rightarrow S)$, which is a usual column insertion, we have
$u_{S'}=u_{S}u_{a'}$. Hence
$$u_{(S,T)}u_a=u_{S}u_{T}u_a=u_{S}u_{a'}u_{T'}=u_{S'}u_{T'}=u_{(a\rightarrow(S,T))}.$$

Next, consider $((S,T)\leftarrow a^\vee)=(S',T')$. If the pair
$(S,(T^\vee \leftarrow a)^\vee)$ satisfies the condition (2) in
Definition \ref{Bmunu}, then $(S',T')=(S,(T^\vee \leftarrow
a)^\vee)$, which implies that $u_{S}=u_{S'}$ and
$u_{T}u_{a^\vee}=u_{T'}$. Hence,
$u_{(S,T)}u_{a^\vee}=u_{((S,T)\leftarrow a^\vee)}$.

Suppose that there exists $j$ such that $a_j=k$ is bumped out of the
$(j-1)$-st row in the row insertion of $a$ into $T^\vee$ and the
insertion of $a_j$ into the $j$-th row violates the condition (2) in
Definition \ref{Bmunu}.

Let $T''=(T^\vee \leftarrow a)^\vee$. Suppose that  $w_{\rm
col}(S)=\widetilde{w}^+w^+$ and $w_{\rm
col}(T'')=w^-\widetilde{w}^-$, where $w^+=w^+_1\ldots w^+_p$ is the
subword  corresponding to the leftmost column of $S$ and
$w^-=w^-_q\ldots w^-_1$ is the subword corresponding to the
rightmost column of $T''$ reading from top to bottom. Note that
$w^-_j=k^\vee$. Suppose that $w^+_i=k$.  By Lemma \ref{reduction},
we have
$$u_{w^+}u_{w^-}=t\,u_{w^+_1\ldots\widehat{{w^+_i}}\ldots
w^+_p}u_{w^-_q \ldots\widehat{w^-_j}\ldots w^-_1}.$$ Note that
$w_{\rm col}(T')=w^-_q \ldots\widehat{w^-_j}\ldots
w^-_1{\widetilde{w}^-}$. Recalling that $S'$ is obtained by playing
the jeu de taquin after removing $k$ in the first column of $S$, it
follows that
$u_{S'}=u_{\widetilde{w}^+}u_{w^+_1\ldots\widehat{{w^+_i}}\ldots
w^+_p}$. Therefore, {\allowdisplaybreaks
\begin{align*}
u_{(S,T)}u_{a^\vee}&=u_{S}u_{T}u_{a^\vee}=u_{S}u_{T''} \\
&=u_{\widetilde{w}^+}u_{w^+}u_{w^-}u_{\widetilde{w}^-}\\
&=t\, u_{\widetilde{w}^+}u_{w^+_1\ldots\widehat{{w^+_i}}\ldots
w^+_p}u_{w^-_q \ldots\widehat{w^-_j}\ldots
w^-_1}u_{\widetilde{w}^-}\\
&=t\,u_{S'}u_{T'}=t\,u_{((S,T)\leftarrow a^\vee)}.
\end{align*}
\qed

Now, we obtain the following immediately.
\begin{prop}\label{plactic rel} For $w=w_1\ldots w_r\in \W$, there exists  $(S,T)\in \B_{\mu,\nu}$ such that
$u_w=t^\epsilon u_{(S,T)}$ where $\epsilon=r-|\mu|-|\nu|$.
\end{prop}

\begin{cor}\label{plactic basis} The set
$\left\{ \,u_{(S,T)}\,\big\vert\,(S,T)\in\B_{\mu,\nu}, \
\mu,\nu\in\cP \,\right\}$ spans $\mathscr{U}(t)$ over
$\mathbb{Q}[t,t^{-1}]$.
\end{cor}

The uniqueness of $(S,T)$ in Proposition \ref{plactic rel} and the
linear independence of the spanning set in Corollary \ref{plactic
basis} will be proved in Section \ref{crystal equivalence}.

\subsection{Crystal equivalence}\label{crystal equivalence} Let $P=P_{\mathbb{N}}=\bigoplus_{i\geq 1}\mathbb{Z}\epsilon_i$ be the
weight lattice of $\gl_{>0}$ and
$\{\,\alpha_i=\epsilon_i-\epsilon_{i+1}\,|\,i\geq 1\,\}$  the set of
simple roots of $\gl_{>0}$. A (regular) $\gl_{>0}$-crystal is a set
$B$ together with the maps ${\rm wt} : B \rightarrow P$,
$\varepsilon_i, \varphi_i: B \rightarrow \mathbb{Z}_{\geq 0}$ and
$\te_i, \tf_i: B \rightarrow B\cup\{{\bf 0}\}$ ($i\geq 1$) such that
for $b\in B$ and $i\geq 1$
\begin{itemize}
\item[(1)] $\varphi_i(b) =\langle {\rm wt}(b),h_i \rangle +
\varepsilon_i(b),$

\item[(2)] $\varepsilon_i(b)={\rm max}\{\,k\,|\,\te_i^kb\neq
{\bf 0}\,\}$ and $\varphi_i(b)={\rm max}\{\,k\,|\,\tf_i^kb\neq {\bf
0}\,\}$

\item[(3)]  ${\rm wt}({\te_i}b)={\rm wt}(b)+\alpha_i$ if
$\te_i b \neq {\bf 0}$, and ${\rm wt}({\tf_i}b)={\rm
wt}(b)-\alpha_i$ if $\tf_i b \neq {\bf 0}$,

\item[(4)] $\tf_i b = b'$ if and only if $b = \te_i
b'$ for $b, b' \in B$,
\end{itemize}
where ${\bf 0}$ is a formal symbol (cf. \cite{Kas94}). Note that $B$
is equipped with a colored oriented graph structure, where
$b\stackrel{i}{\rightarrow}b'$ if and only if $b'=\tf_{i}b$ for
$b,b'\in B$ and $i\geq 1$. The dual crystal $B^\vee$ of $B$ is
defined to be the set $\{\,b^\vee\,|\,b\in B\,\}$ with ${\rm
wt}(b^\vee)=-{\rm wt}(b)$,
$\widetilde{e}_i(b^\vee)=\left(\widetilde{f}_i b \right)^\vee$ and
$\widetilde{f}_i(b^\vee)=\left(\widetilde{e}_i b \right)^\vee$ for
$b\in B$ and $i\geq 1$. We assume that ${\bf 0}^\vee={\bf 0}$. Note
that $\mathbb{N}$ is naturally equipped with a $\gl_{>0}$-crystal
structure;
$$ 1 \stackrel{1}{\longrightarrow} 2
\stackrel{2}{\longrightarrow}3\stackrel{3}{\longrightarrow} \cdots$$
with ${\rm wt}(i)=\epsilon_i$ ($\i\geq 1$), while $\mathbb{N}^\vee$
is its dual.

For crystals $B_1$ and $B_2$, a tensor product $B_1\otimes B_2$ is
defined to be the set $\{\,b_1\otimes b_2\,|\,b_i\in B_i\,\,
(i=1,2)\,\}$ with ${\rm wt}(b_1\otimes b_2)={\rm wt}(b_1)+{\rm
wt}(b_2)$, and
\begin{equation*}
{\te}_i(b_1\otimes b_2)=
\begin{cases}
{\te}_i b_1 \otimes b_2, & \text{if $\varphi_i(b_1)\geq \varepsilon_i(b_2)$}, \\
b_1\otimes {\te}_i b_2, & \text{if
$\varphi_i(b_1)<\varepsilon_i(b_2)$},
\end{cases}
\end{equation*}
\begin{equation*}
{\tf}_i(b_1\otimes b_2)=
\begin{cases}
{\tf}_i b_1 \otimes b_2, & \text{if  $\varphi_i(b_1)>\varepsilon_i(b_2)$}, \\
b_1\otimes {\tf}_i b_2, & \text{if $\varphi_i(b_1)\leq
\varepsilon_i(b_2)$},
\end{cases}
\end{equation*}
for $i\geq 1$ and $b_1\otimes b_2\in B_1\otimes B_2$, where we
assume that ${\bf 0}\otimes b_2=b_1\otimes {\bf 0}={\bf 0}$. For
example, $\W$ is a $\gl_{>0}$-crystal, where each word $w_1\ldots
w_r$ is identified with $w_1\otimes \cdots \otimes w_r$ in a mixed
$r$-tensor product of $\mathbb{N}$ and $\mathbb{N}^\vee$.

For $b_i\in B_i$ ($i=1,2$), we say that $b_1$ is equivalent to
$b_2$, and write $b_1 \equiv b_2$ if ${\rm wt}(b_1)={\rm wt}(b_2)$
and they generate the same $\mathbb{N}$-colored graph with respect
to $\te_i$, $\tf_i$ ($i\geq 1$). We usually call $\equiv$ the
crystal equivalence..

For a skew Young diagram $\lambda/\mu$, $\B_{\lambda/\mu}$ has a
well-defined $\gl_{>0}$-crystal structure such that
$\widetilde{x}_i(S)=S'$ if $\widetilde{x}_i w_{\rm col}(S)\neq {\bf
0}$ ($i\geq 1$, $x=e,f$), where $S'$ is the unique tableau in
$\B_{\lambda/\mu}$ with $w_{\rm col}(S')=\widetilde{x}_i w_{\rm
col}(S)$, and $\widetilde{x}_i(S)={\bf 0}$ otherwise. We regard
$\B_{\lambda/\mu}^\vee$ as the dual of $\B_{\lambda/\mu}$. Moreover,
for $\mu,\nu\in\cP$, $\B_{\mu,\nu}\cup\{{\bf 0}\}\subset
\left(\B_\mu\otimes \B_\nu^\vee\right)\cup\{{\bf 0}\}$ is invariant
under $\te_i$, $\tf_i$ ($i\geq 1$), and hence a $\gl_{>0}$-crystal,
which is connected as a graph \cite[Proposition 3.4]{K09}.

Let $\mathscr{W}$ be an associative $\mathbb{Q}$-algebra generated
by the symbol $[w]$ ($w\in\W$) subject to the relations;
\begin{equation*}\label{crystal plactic}
\begin{split}
[w][w']&=[w w'], \\
[w]&=[w'], \ \ \ \ \text{if $w\equiv w'$},
\end{split}
\end{equation*}
for $w, w'\in\W$. Note that $[\emptyset]=1$ is the unity in
$\mathscr{W}$, where $\emptyset$ is the empty word.

\begin{lem}\label{basis of W} The set
$$\mathscr{B}=
\left\{\,\bigl[w_{\rm col}(S)w_{\rm
col}(T)\bigr]\,\big\vert\,(S,T)\in\B_{\mu,\nu}, \ \mu,\nu\in
\cP\,\right\}$$ is a $\mathbb{Q}$-basis of $\mathscr{W}$.
\end{lem}
\pf For $a\in \mathbb{N}$ and $(S,T)\in\B_{\mu,\nu}$, it is shown in
\cite{K09} that$$\left(a\rightarrow (S,T)\right)\,\equiv
\,(S,T)\otimes a, \ \ \left((S,T)\leftarrow a^\vee\right) \,\equiv\,
(S,T)\otimes a^\vee.$$ This implies that for $w\in\W$,
$[w]=\left[w_{\rm col}(S)w_{\rm col}(T)\right]$ for some
$(S,T)\in\B_{\mu,\nu}$, and hence $\mathscr{W}$ is spanned by
$\mathscr{B}$.

Now, suppose that
\begin{equation}\label{linear independence}
\sum_{i=1}^n c_i \left[w_{\rm col}(S^{(i)})w_{\rm
col}(T^{(i)})\right]=0
\end{equation}
for some $c_i\in\mathbb{Q}$ and $(S^{(i)},T^{(i)})\in
\B_{\mu^{(i)},\nu^{(i)}}$ with $\mu^{(i)},\nu^{(i)}\in\cP$ ($1\leq
i\leq n$). Since $(S,T)\equiv (S',T')$ implies $(S,T)=(S',T')$ for
$(S,T)\in\B_{\mu,\nu}$ and $(S',T')\in\B_{\sigma,\tau}$ \cite[Lemma
5.1]{K09}, we assume that $(S^{(i)},T^{(i)})$'s are mutually
different.

We use induction on $n$ to show that $c_i=0$ for $1\leq i\leq n$. It
is clear when $n=1$. Suppose  that $n\geq 2$. Since no two of
$(S^{(i)},T^{(i)})$'s are mutually crystal equivalent, there exist
$j_1,\ldots,j_r$ such that $\widetilde{x}_{j_1}\cdots
\widetilde{x}_{j_r}(S^{(1)},T^{(1)})={\bf 0}$ but
$\widetilde{x}_{j_1}\cdots \widetilde{x}_{j_r}(S^{(i)},T^{(i)})\neq
{\bf 0}$ for some $2\leq i\leq n$, where $x=e$ or $f$ for each $j_k$
($1\leq k\leq r$).

Note that $\widetilde{x}_i$ ($x=e,f$, $i\geq 1$) acts on
$\mathscr{W}$ by $\widetilde{x}_i[w]=[\widetilde{x}_iw]$, where we
assume that $[{\bf 0}]=0$. Hence by applying
$X=\widetilde{x}_{j_1}\cdots \widetilde{x}_{j_r}$ to (\ref{linear
independence}), we get
\begin{equation*}\label{linear independence'}
\sum_{i=2}^n c_i \left[Xw_{\rm col}(S^{(i)})w_{\rm
col}(T^{(i)})\right]=\sum_{i=2}^n c_i \left[w_{\rm
col}(\overline{S}^{(i)})w_{\rm col}(\overline{T}^{(i)})\right]=0
\end{equation*}
for some $\bigl[w_{\rm col}(\overline{S}^{(i)})w_{\rm
col}(\overline{T}^{(i)})\bigr]\in \mathscr{B}$.

Here, we assume that $c_i=0$ if $X\left(w_{\rm col}(S^{(i)})w_{\rm
col}(T^{(i)})\right)={\bf 0}$. By induction hypothesis, we have
$c_2=\cdots=c_n=0$, and hence $c_1=0$. Therefore, $\mathscr{B}$ is a
$\mathbb{Q}$-basis of $\mathscr{W}$. \qed

\begin{thm}\label{UisoK} Let $\mathscr{U}(1)$ be the $\mathbb{Q}$-algebra obtained from
$\mathscr{U}(t)$ by specializing $t=1$. Then as a
$\mathbb{Q}$-algebra, we have
$$\mathscr{U}(1) \simeq \mathscr{W},$$
where $u_a$ is mapped to $[a]$ for $a\in
\mathbb{N}\cup\mathbb{N}^\vee$.
\end{thm}
\pf  It is straightforward to see that the relations in
(\ref{plactic rel}) are satisfied in $\mathscr{W}$ under the
correspondence $u_a\mapsto [a]$. Hence there exists a
$\mathbb{Q}$-algebra  homomorphism  $\psi :
\mathscr{U}(1)\rightarrow \mathscr{W}$ sending $u_a$ to $[a]$. Since
$\left\{ \,u_{(S,T)}\,\big\vert\,(S,T)\in\B_{\mu,\nu}, \
\mu,\nu\in\cP \,\right\}$ spans $\mathscr{U}(1)$ and
$\psi(u_{(S,T)})=\left[w_{\rm col}(S)w_{\rm col}(T)\right]$, it
follows from Lemma \ref{basis of W} that $\psi$ is an isomorphism.
\qed

\begin{cor}\label{basis of U}
The set
\begin{equation*}
\begin{split}
&\left\{ \,u_{(S,T)}\,\big\vert\,(S,T)\in\B_{\mu,\nu}, \
\mu,\nu\in\cP \,\right\}
\end{split}
\end{equation*}
is a $\mathbb{Q}[t,t^{-1}]$-basis of $\mathscr{U}(t)$.
\end{cor}
\pf Since $\left\{ \,u_{(S,T)}\,\big\vert\,(S,T)\in\B_{\mu,\nu}, \
\mu,\nu\in\cP \,\right\}\subset \mathscr{U}(1)$ is mapped to
$\mathscr{B}$ by Theorem \ref{UisoK}, it is a $\mathbb{Q}$-basis of
$\mathscr{U}(1)$. Now  there is a well-defined $\mathbb{Q}$-algebra
homomorphism $\psi' : \mathscr{U}(t)\rightarrow \mathscr{U}(1)$ such
that $\psi'(u_a)=u_a$ and $\psi'(t)=1$. Then it is not difficult to
check that $\left\{ \,u_{(S,T)}\,\big\vert\,(S,T)\in\B_{\mu,\nu}, \
\mu,\nu\in\cP \,\right\}\subset \mathscr{U}(t)$ is linearly
independent over $\mathbb{Q}[t,t^{-1}]$ and hence a
$\mathbb{Q}[t,t^{-1}]$-basis of $\mathscr{U}(t)$. \qed

\begin{cor}
For $w\in\W$, there exists unique $(S,T)\in \B_{\mu,\nu}$  and
$\epsilon\in\Z_{\geq 0}$ such that $u_w=t^\epsilon u_{(S,T)}$.
\end{cor}

\subsection{Non-commutative Schur functions}
Let $\widehat{\mathscr{U}(t)}=\bigoplus_{n\geq
0}\widehat{\mathscr{U}(t)}_n$, where $\widehat{\mathscr{U}(t)}_n$ is
the completion of
$\bigoplus_{|\mu|+|\nu|=n}\bigoplus_{(S,T)\in\B_{\mu,\nu}}\mathbb{Q}[t,t^{-1}]u_{(S,T)}$.
For a skew Young diagram $\lambda/\mu$, let
\begin{equation*}
\begin{split}
s_{\lambda/\mu}({\bf u})&=\sum_{S\in \B_{\lambda/\mu}}u_S, \ \
s_{\lambda/\mu}^\vee({\bf u})=\sum_{S\in \B_{\lambda/\mu}^\vee}u_S\
\ \in \widehat{\mathscr{U}(t)},
\end{split}
\end{equation*}
which are plactic skew Schur functions in $u_i$'s and
$u_{i^\vee}$'s, respectively.

Let $\Lambda(t)$ be the algebra of symmetric functions in ${\bf
x}={\bf x}_{\mathbb{N}}$ over $\mathbb{Q}[t,t^{-1}]$. Then
$\{\,s_{(k)}({\bf u})\,|\,k\geq 1\,\}$ (resp. $\{\,s_{(k)}^\vee({\bf
u})\,|\,k\geq 1\,\}$) generates the subalgebra $\mathscr{S}(t)_\pm$
of $\widehat{\mathscr{U}(t)}$ isomorphic to $\Lambda(t)$ \cite{LS},
where $s_{(k)}({\bf u})$ (resp. $s_{(k)}^\vee({\bf u})$) corresponds
to the $k$-th complete symmetric function $h_k({\bf x})=s_{(k)}({\bf
x})$, and $\{\,s_\lambda({\bf u})\,|\,\lambda\in\cP\,\}$ (resp.
$\{\,s_\lambda^\vee({\bf u})\,|\,\lambda\in\cP\,\}$) is a
$\mathbb{Q}[t,t^{-1}]$-basis of $\mathscr{S}(t)_+$ (resp.
$\mathscr{S}(t)_-$).

We define
\begin{equation*}
s_{\mu,\nu}({\bf u})=\sum_{(S,T)\in\B_{\mu,\nu}}u_{(S,T)}
\end{equation*}
for $\mu,\nu\in\cP$ and let
\begin{equation*}
\mathscr{S}(t)=\sum_{\mu,\nu\in\cP}\mathbb{Q}[t,t^{-1}]s_{\mu,\nu}({\bf
u})\subset \widehat{\mathscr{U}(t)}.
\end{equation*}

\begin{lem}\label{plactic Schur rel} For $\mu,\nu\in\cP$, we have
\begin{equation*}
\begin{split}
&s_\mu({\bf u})s_\nu^\vee({\bf u})=\sum_{\lambda\subset
\mu,\nu}t^{|\lambda|}s_{\nu/\lambda}^\vee({\bf
u})s_{\mu/\lambda}({\bf u})
=\sum_{\lambda,\sigma,\tau}t^{|\lambda|}c^{\mu}_{\lambda\,\sigma}c^{\nu}_{\lambda\,\tau}s_{\sigma,\tau}({\bf
u}).
\end{split}
\end{equation*}
\end{lem}
\pf By (\ref{bijection-1}) and Lemma \ref{plactic lemma} (1), we
have $s_{\mu,\nu}({\bf u})=s_{\nu}^\vee({\bf u})s_\mu({\bf u})$. The
identity follows from (\ref{extLR}) and Lemma \ref{plactic lemma}
(2). \qed

\begin{prop}
$\mathscr{S}(t)$ is a $\mathbb{Q}[t,t^{-1}]$-algebra with a basis
$\{\,s_{\mu,\nu}({\bf u})\,|\,\mu,\nu\in\cP\,\}$, where
\begin{equation*}
s_{\mu,\nu}({\bf u}) s_{\sigma,\tau}({\bf u})=
\sum_{\zeta,\eta}\left(\sum_{\alpha,\beta,\gamma}t^{|\gamma|}
c^{\zeta}_{\sigma\,\alpha}c^{\mu}_{\alpha\beta}c^{\tau}_{\beta\,\gamma}c^{\eta}_{\gamma\nu}\right)
s_{\zeta,\eta}({\bf u}).
\end{equation*}
for $\mu,\nu,\sigma,\tau\in\cP$.
\end{prop}
\pf In fact, $\{\,s_{\mu,\nu}({\bf u})\,|\,\mu,\nu\in\cP\,\}$ is
linearly independent  by Lemma \ref{basis of U}, and hence a basis
of $\mathscr{S}(t)$. Combining Lemma \ref{plactic Schur rel} with
the usual Littlewood-Richardson rule (\ref{LR}) for $s_\mu({\bf
u})$'s and $s_\nu^\vee({\bf u})$'s, we obtain the above identity.
Since the sum on the right hand side is finite, $\mathscr{S}(t)$ has
a well-defined multiplication and hence is a
$\mathbb{Q}[t,t^{-1}]$-algebra. \qed

\subsection{Heisenberg algebra}
Let $\mathscr{H}(t)$ be an associative
$\mathbb{Q}[t,t^{-1}]$-algebra with unity generated by $B_n$
($n\in\Z\setminus\{0\}$) subject to the relations;
\begin{equation*}
B_kB_l-B_lB_k = k\,t^k\delta_{k+l,0}.
\end{equation*}
For $k\geq 1$, let
$p_k({\bf u})\in \mathscr{S}(t)_+$ (resp. $p^\vee_k({\bf u})\in
\mathscr{S}(t)_-$) correspond to the $k$-th power sum symmetric
function $p_k({\bf x})\in\Lambda(t)$.

\begin{prop}
As a $\mathbb{Q}[t,t^{-1}]$-algebra, we have $$\mathscr{S}(t)\simeq
\mathscr{H}(t),$$ where $p_k({\bf u})$ and $p^\vee_k({\bf u})$ are
mapped to $B_k$ and $B_{-k}$ for $k\geq 1$, respectively.
\end{prop}
\pf  By Lemma \ref{plactic Schur rel}, we have
\begin{equation}\label{defining rel for K}
h_s({\bf u}) h^\vee_r({\bf u})=\sum_{i=0}^m t^i\,h^\vee_{r-i}({\bf
u})h_{s-i}({\bf u})
\end{equation}
for $r,s\geq 1$, where $m=\min\{r,s\}$. Here we assume that
$h_0({\bf u})=h^\vee_0({\bf u})=1$. We may view $\mathscr{S}(t)$ as
an algebra generated by $\{\,h^\vee_s({\bf u}), h_r({\bf
u})\,|\,r,s\geq 1\,\}$ with the defining relations (\ref{defining
rel for K}). Since
\begin{equation*}
h_r({\bf u})=\sum_{|\lambda|=r}\frac{1}{z_\lambda}p_{\lambda}({\bf
u}),
\end{equation*}
where $z_\lambda=\prod_{i\geq 1}i^{m_i(\lambda)}m_i(\lambda)!$ and
$m_i(\lambda)=\left|\{\,k\,|\,\lambda_k=i\,\}\right|$,  we obtain
$$p_k({\bf u})p^\vee_l({\bf u})-p^\vee_l({\bf u})p_k({\bf u}) =
k\,t^k\delta_{k,l}$$ by using the same argument as in
\cite[Corollary 8]{L}. This implies that the map $\psi :
\mathscr{H}(t)\rightarrow \mathscr{S}(t)$ sending $B_{-k}$ (resp.
$B_k$) to $p^\vee_k({\bf u})$ (resp. $p_k({\bf u})$) for $k\geq 1$
is a well-defined isomorphism.\qed

\begin{rem}{\rm
Regarding $\mathscr{S}(0)$ and $\mathscr{S}(1)$ as
$\mathbb{Q}$-algebras generated by $h_k({\bf u})$ and $h^\vee_k({\bf
u})$ ($k\geq 1$), we have $\mathscr{S}(0)\simeq \Lambda\otimes
\Lambda$, and $\mathscr{S}(1)\simeq \bigl\langle\,
\frac{\partial}{\partial p_k}, p_k\,|\,k\geq 1\,\bigr\rangle\subset
{\rm End}_{\mathbb{Q}}(\Lambda)$, where $\Lambda$ is the algebra of
symmetric functions in ${\bf x}$ over $\mathbb{Q}$ and $p_k$ is the
operator on $\Lambda$ induced from the multiplication by $p_k({\bf
x})$. Therefore, we may view $\mathscr{S}(t)$ as an algebra
interpolating the algebra of double symmetric functions and the Weyl
algebra of infinite rank.}
\end{rem}

\section{Knuth correspondence and Cauchy identity}

\subsection{Main result}
Let $\A$ and $\B$ be linearly ordered $\Z_2$-graded sets. For $A\in \M_{\A,\mathbb{N}}$ (or $\M_{\A,\mathbb{N}^\vee}$), we put $u_A=u_{\bj}$ if $A=A(\bi,\bj)$. Now we are
in a position to state and prove our main theorem.
\begin{thm}\label{main result} There exists a bijection
$$\M_{\A,\mathbb{N}}\times \M_{\B,\mathbb{N}^\vee}
\longrightarrow \M_{\B,\mathbb{N}^\vee}\times \M_{\A,\mathbb{N}}
\times \M_{\A,\B}$$ sending $(X,Y)$ to $(Y',X',Z)$ such that
\begin{itemize}
\item[(1)] ${\rm wt}_{\A}(X)={\rm wt}_{\A}(X')+{\rm wt}_{\A}(Z)$,  ${\rm wt}_{\B}(Y)={\rm wt}_{\B}(Y')+{\rm wt}_{\B}(Z)$,

\item[(2)] $u_{X}u_Y=t^{|Z|}u_{Y'}u_{X'}$ where $Z=(z_{ij})$ and $|Z|=\sum_{i,j}z_{ij}$.
\end{itemize}

\end{thm}
\pf It is obtained by composing the following
bijections;
{\allowdisplaybreaks
\begin{align*}
&\M_{\A,\mathbb{N}}\times \M_{\B,\mathbb{N}^\vee}\\
&{\longrightarrow} \bigsqcup_{\mu\in \cP_\A, \nu\in \cP_\B}\B_{\mu}
\times SST_\A(\mu)\times \B_{\nu}^\vee \times SST_\B(\nu)\ \ \ \ \ \
\ \ \ \ \ \ \ \ \ \ \text{by (\ref{Knuth-1}) and (\ref{Knuth-2})} \\
&{\longrightarrow} \bigsqcup_{\mu\in \cP_\A, \nu\in
\cP_\B}SST_\A(\mu)\times  SST_\B(\nu)\times \B_{\mu} \times
\B_{\nu}^\vee\\
&{\longrightarrow} \bigsqcup_{\mu\in \cP_\A, \nu\in
\cP_\B}SST_\A(\mu)\times SST_\B(\nu)\times
\left(\bigsqcup_{\lambda\subset \mu,\nu}\cB_{\nu/\lambda}^\vee\times
\cB_{\mu/\lambda}\right)\ \ \ \ \ \ \ \  \text{by (\ref{extLR})}\\
&{\longrightarrow} \bigsqcup_{\mu\in \cP_\A, \nu\in
\cP_\B}\bigsqcup_{\lambda\subset \mu,\nu} SST_\B(\nu)\times
\cB_{\nu/\lambda}^\vee\times
SST_\A(\mu)\times \cB_{\mu/\lambda}\\
&{\longrightarrow} \bigsqcup_{\lambda\in \cP_\A\cap \cP_\B}
\M_{\B,\mathbb{N}^\vee} \times SST_\B(\lambda)\times
\M_{\A,\mathbb{N}} \times
SST_\A(\lambda) \ \ \ \ \ \ \ \ \ \ \ \ \ \ \ \ \ \ \ \ \text{by (\ref{Knuth-3})}\\
&{\longrightarrow}\ \ \M_{\B,\mathbb{N}^\vee}\times
\M_{\A,\mathbb{N}} \times \left(\bigsqcup_{\lambda\in \cP_\A\cap
\cP_\B}
SST_\B(\lambda)\times SST_\A(\lambda)\right)\\
&{\longrightarrow}\ \ \M_{\B,\mathbb{N}^\vee}\times
\M_{\A,\mathbb{N} } \times \M_{\A,\B} \ \ \ \ \ \ \ \ \ \ \ \ \ \ \ \ \ \ \  \ \ \ \ \ \ \ \ \ \ \ \ \ \ \ \ \ \ \ \ \ \ \ \ \ \ \ \ \ \ \ \ \text{by (\ref{Knuth-1})}.\\
\end{align*}}
\qed\vskip 2mm

Now, let us consider the non-commutative character identity
associated with Theorem \ref{main result}. We first define the
plactic Cauchy products
\begin{equation*}
\mathscr{Q}({\bf
x}_\A)=\overrightarrow{\prod_{a\in\A}}\mathscr{Q}(x_a), \ \ \
\mathscr{P}({\bf
x}_\B)=\overrightarrow{\prod_{b\in\B}}\mathscr{P}(x_b),
\end{equation*}
where the products are given with respect to the linear ordering on
$\A$ or $\B$ so that smaller terms are to the left, and
{\allowdisplaybreaks \begin{align*}
\mathscr{P}(x_b)&=
\begin{cases}
\dfrac{1}{\cdots(1-u_{2^\vee}x_b)(1-u_{1^\vee}x_b)}, & \text{if
$|b|=0$}, \\
{(1+u_{1^\vee}x_b)(1+u_{2^\vee}x_b)\cdots}, & \text{if $|b|=1$},
\end{cases}\\
\mathscr{Q}(x_a)&=\begin{cases}
\dfrac{1}{\cdots(1-u_{2}x_a)(1-u_{1}x_a)}, & \text{if
$|a|=0$}, \\
{(1+u_{1}x_a)(1+u_{2}x_a)\cdots}, & \text{if $|a|=1$}.
\end{cases}
\end{align*}}
We assume that ${\bf x}_\A$ and ${\bf x}_\B$ commute with ${\bf u}$.
Note that
\begin{equation*}\label{Cauchy-1}
\mathscr{Q}({\bf x}_\A)=\sum_{\lambda\in\cP_\A}s_\lambda({\bf
u})s_\lambda({\bf x}_\A), \ \ \mathscr{P}({\bf
x}_\B)=\sum_{\lambda\in\cP_\B}s_\lambda^\vee({\bf u})s_\lambda({\bf
x}_\B),
\end{equation*}
by (\ref{Knuth-1}) and (\ref{Knuth-2}). Since the bijections in the
proof of Theorem \ref{main result} preserves the plactic relations
(\ref{plactic}), ${\rm wt}_\A$ and ${\rm wt}_\B$, we obtain the
following identity.

\begin{cor}\label{Cauchy identity}
$$\mathscr{Q}({\bf x}_\A)\mathscr{P}({\bf x}_\B)=\mathscr{P}({\bf x}_\B)\mathscr{Q}({\bf x}_\A)
\frac{\prod_{|a|\neq |b|}(1+tx_ax_b)}{\prod_{|a|=
|b|}(1-tx_ax_b)}.$$
\end{cor}

\subsection{Cauchy identity for Schur operators}
For $i\in\mathbb{N}$, we define operators $\overline{u}_i,
\overline{u}_{i^\vee}\in {\rm
End}_{\mathbb{Q}[t,t^{-1}]}(\Lambda(t))$ by
\begin{equation*}
\begin{split}
\overline{u}_{i^\vee}(s_{\mu}({\bf x}))&=
\begin{cases}
s_{\mu\cup\{(i,\mu_i+1)\}}({\bf x}), & \text{if
$\mu\cup\{(i,\mu_i+1)\}\in\cP$},\\
0, & \text{if $\mu\cup\{(i,\mu_i+1)\}\not\in\cP$},
\end{cases}\\
\overline{u}_{i}(s_{\mu}({\bf x}))&=
\begin{cases}
t\,s_{\mu\setminus\{(i,\mu_i)\}}({\bf x}), & \text{if
$\mu\setminus\{(i,\mu_i)\}\in\cP$},\\
0, & \text{if $\mu\setminus\{(i,\mu_i)\}\not\in\cP$}.
\end{cases}
\end{split}
\end{equation*}
These operators are called {Schur operators} \cite{Fomin95}. Let
$\overline{\mathscr{U}}(t)$ be the subalgebra of ${\rm
End}_{\mathbb{Q}[t,t^{-1}]}(\Lambda(t))$ generated by
$\overline{u}_i, \overline{u}_{i^\vee}$ ($i\geq 1$). It is easy to
see that there exists a surjective $\mathbb{Q}[t,t^{-1}]$-algebra
homomorphism $\psi : \mathscr{U}(t) \rightarrow
\overline{\mathscr{U}}(t)$ such that $\psi(u_i)=\overline{u}_{i}$
and $\psi(u_{i^\vee})=\overline{u}_{i^\vee}$ for $i\in\mathbb{N}$.

For $\lambda\in \cP$, let
\begin{equation*}
s_\lambda(\overline{\bf u})=\sum_{S\in \B_\lambda}\overline{u}_S, \
\ s^\vee_\lambda(\overline{\bf u})=\sum_{S\in
\B_\lambda^\vee}\overline{u}_S,
\end{equation*}
where $\overline{u}_S=\psi(u_S)$ for $S\in \B_\lambda$ or
$\B_\lambda^\vee$. For $\lambda,\mu\in\cP$, we have
\begin{equation*}
s^\vee_\mu(\overline{\bf u})(s_{\lambda}({\bf x}))=s_\lambda({\bf
x})s_{\mu}({\bf x}), \ \ s_\mu(\overline{\bf u})(s_{\lambda}({\bf
x}))=t^{|\mu|}s_{\lambda/\mu}(\bf x)
\end{equation*}
(see \cite{Fomin95}). We also have
\begin{equation*}
\begin{split}
&\overline{\mathscr{Q}}({\bf
x}_\A)=\overrightarrow{\prod_{a\in\A}}\overline{\mathscr{Q}}(x_a)
=\sum_{\lambda\in\cP_\A}s_\lambda(\overline{\bf u}) s_\lambda({\bf
x}_\A), \ \\
&\overline{\mathscr{P}}({\bf
x}_\B)=\overrightarrow{\prod_{b\in\B}}\overline{\mathscr{P}}(x_b)
=\sum_{\lambda\in\cP_\B}s^\vee_\lambda(\overline{\bf u})
s_\lambda({\bf x}_\B),
\end{split}
\end{equation*}
where $\overline{\mathscr{P}}(x_a)$ and
$\overline{\mathscr{Q}}(x_b)$ are obtained from ${\mathscr{P}}(x_a)$
and ${\mathscr{Q}}(x_b)$ by replacing $u_i, u_{i^\vee}$ with $\overline{u}_{i},\overline{u}_{i^\vee}$,
respectively. Therefore, the products $\overline{\mathscr{Q}}({\bf
x}_\A)\overline{\mathscr{P}}({\bf x}_\B)$ and
$\overline{\mathscr{P}}({\bf x}_\B)\overline{\mathscr{Q}}({\bf
x}_\A)$ are well defined, and the identity in Corollary \ref{Cauchy
identity} gives the following, which recovers the generalized Cauchy
identity for Schur operators \cite{Fomin95} when $t=1$;
$$\overline{\mathscr{Q}}({\bf x}_\A)\overline{\mathscr{P}}({\bf x}_\B)
=\overline{\mathscr{P}}({\bf x}_\B)\overline{\mathscr{Q}}({\bf
x}_\A) \frac{\prod_{|a|\neq |b|}(1+tx_ax_b)}{\prod_{|a|=
|b|}(1-tx_ax_b)}.$$

\subsection{Knuth correspondence for skew tableaux}\label{skewRSK}
Fix $\alpha,\beta \in\cP$. For $w=w_1\ldots w_r\in \W$, let us say
that $w$ is a Littlewood-Richardson (simply LR) word  of shape
$(\alpha,\beta)$ if there exist $\alpha^{(k)}\in\cP$ for $1\leq
k\leq r$ such that $(w_1\ldots w_k \rightarrow
H_{\alpha})=H_{\alpha^{(k)}}$ and $\alpha^{(r)}=\beta$. Note that
for $0\leq k\leq r-1$, $|\alpha^{(k+1)}|=|\alpha^{(k)}|+ 1$ if
$w_k\in\mathbb{N}$ and $|\alpha^{(k+1)}|=|\alpha^{(k)}|- 1$ if
$w_k\in\mathbb{N}^\vee$ (we assume that $\alpha^{(0)}=\alpha$), and
by definition $w_1\ldots w_k$ is also a LR word of shape
$(\alpha,\alpha^{(k)})$ for $1\leq k\leq r$.

\begin{lem}\label{LR word} Let $\alpha, \beta\in\cP$ be given. For $w\in \W$,
$w$ is a LR word of shape $(\alpha,\beta)$ if and only if
$H_\alpha\otimes w \equiv H_\beta$. In particular, $w'$ is a LR word
of shape $(\alpha,\beta)$ for all $w'\in\W$ with $w'\equiv w$.
\end{lem}
\pf Suppose that $w$ is a LR word of shape $(\alpha,\beta)$. Since
$(w_1\ldots w_k \rightarrow H_{\alpha})=H_\alpha \otimes w_1\ldots
w_k$ for $1\leq k\leq r$, it is clear that $H_\alpha\otimes w \equiv
H_\beta$.

Conversely, suppose that  $H_\alpha\otimes w \equiv H_\beta$. If
$(w_1\ldots w_k \rightarrow H_{\alpha})\in \B_{\mu,\nu}$ for some
$k$ and $\mu,\nu\in\cP$ with $\nu\neq \emptyset$, then $(w
\rightarrow H_\alpha) \in\B_{\sigma,\tau}$ for some
$\sigma,\tau\in\cP$ with $\tau\neq \emptyset$, which contradicts the
fact that $H_\alpha\otimes w\equiv H_\beta$. Hence there exist
$\alpha^{(k)}\in\cP$ for $1\leq k\leq r$ such that $(w_1\ldots w_k
\rightarrow H_{\alpha})\in \B_{\alpha^{(k)}}$ and
$\alpha^{(r)}=\beta$.

Now suppose that $(w_1\ldots w_k \rightarrow H_\alpha)\neq
H_{\alpha^{(k)}}$ for some $k\geq 1$, which is equivalent to saying
that $\te_i((w_1\ldots w_k \rightarrow H_\alpha))\neq {\bf 0}$ for
some $i\geq 1$. Then $\te_i H_\beta\equiv \te_i(w\rightarrow
H_\alpha)\equiv (\te_i(w_1\ldots w_k \rightarrow H_\alpha))\otimes
w_{k+1}\ldots w_r  \neq {\bf 0}$, which is also a contradiction.
This completes the proof. \qed\vskip 2mm

For $\lambda,\mu\in\cP$ with $|\lambda|=|\alpha|+|\mu|$, we have  by
(\ref{LR})
\begin{equation}\label{LRword-1}
\{\,S\in\B_\mu\,|\,w_{\rm col}(S) \text{ is a LR word of shape
$(\alpha,\lambda)$} \,\}\ \stackrel{1-1}{\longleftrightarrow} \ {\bf
LR}^{\lambda}_{\alpha \mu}.
\end{equation}
For $\lambda,\nu\in\cP$ with $|\lambda|=|\beta|+|\nu|$, we have by
(\ref{bijection-2})
\begin{equation}\label{LRword-2}
\{\,S\in\B_\nu^\vee\,|\,w_{\rm col}(S) \text{ is a LR word of shape
$(\lambda,\beta)$} \,\}\ \stackrel{1-1}{\longleftrightarrow}\ {\bf
LR}^{\lambda}_{\beta \nu}.
\end{equation}
Let $\left(\B_{\mu}\times \B_\nu^\vee\right)_{(\alpha,\beta)}$ be
the set of $(S,T)\in \B_{\mu}\times \B_\nu^\vee$ such that $w_{\rm
col}(S)w_{\rm col}(T)$ is a LR word of shape $(\alpha,\beta)$.
Combining (\ref{LRword-1}), (\ref{LRword-2}) and Lemma \ref{LR
word}, we have
\begin{equation}\label{LRword-3}
\left(\B_{\mu}\times \B_\nu^\vee\right)_{(\alpha,\beta)}\
\stackrel{1-1}{\longleftrightarrow}\ \bigsqcup_{\lambda}{\bf
LR}^{\lambda}_{\alpha \mu}\times{\bf LR}^{\lambda}_{\beta \nu}.
\end{equation}
Similarly, for $\sigma, \tau\in\cP$, let $\left(\B_{\tau}^\vee\times
\B_\sigma\right)_{(\alpha,\beta)}$ be the set of $(S,T)\in
\B_{\tau}^\vee\times \B_\sigma$ such that $w_{\rm col}(S)w_{\rm
col}(T)$ is a LR word of shape $(\alpha,\beta)$. As in
(\ref{LRword-3}), we have a bijection
\begin{equation}\label{LRword-4}
\left(\B_{\tau}^\vee\times \B_\sigma\right)_{(\alpha,\beta)}\
\stackrel{1-1}{\longleftrightarrow}\ \bigsqcup_{\lambda}{\bf
LR}^{\alpha}_{\lambda \tau}\times{\bf LR}^{\beta}_{\lambda \sigma}.
\end{equation}

\begin{cor} Let $\alpha,\beta,\mu,\nu \in\cP$ be given. The  bijection
{\rm (\ref{bijection-2})} when restricted to $\left(\B_{\mu}\times
\B_\nu^\vee\right)_{(\alpha,\beta)}$ gives the following bijection
\begin{equation*}
\bigsqcup_{\lambda}{\bf LR}^{\lambda}_{\alpha \mu}\times{\bf
LR}^{\lambda}_{\beta \nu}\ \ {\longrightarrow}
\bigsqcup_{\eta,\zeta, \sigma,\tau}{\bf LR}^{\alpha}_{\eta
\tau}\times{\bf LR}^{\beta}_{\eta \sigma}\times{\bf LR}^\mu_{\sigma
\zeta}\times {\bf LR}^{\nu}_{\tau \zeta}.
\end{equation*}
\end{cor}
\pf Since the bijection (\ref{bijection-2}) preserves the plactic
relations or the crystal equivalence, we have
\begin{equation*}
\left(\cB_\mu\times \cB_\nu^\vee\right)_{(\alpha,\beta)}
{\longrightarrow}
\bigsqcup_{\zeta,\sigma,\tau}\left(\B_{\tau}^\vee\times
\B_\sigma\right)_{(\alpha,\beta)}\times{\bf LR}^\mu_{\sigma
\zeta}\times {\bf LR}^{\nu}_{\tau \zeta}
\end{equation*}
Hence, it follows from (\ref{LRword-3}) and (\ref{LRword-4}).
\qed\vskip 2mm

Let $\left(\M_{\A,\mathbb{N}}\times
\M_{\B,\mathbb{N}^\vee}\right)_{(\alpha,\beta)}$  be the set of
$(A,A')$ such that $\bj \cdot \bj'\in\W$ is a LR word of shape
$(\alpha,\beta)$, where $A=A(\bi,\bj)$ and $A'=A(\bi',\bj')$, and
let $\left( \M_{\B,\mathbb{N}^\vee}\times \M_{\A,\mathbb{N}}
\right)_{(\alpha,\beta)}$ be defined in the same way.

Now, we recover the Knuth type correspondence for skew tableaux by
Sagan and Stanley \cite{SS} as a restriction of Theorem \ref{main
result} to LR words of shape $(\alpha,\beta)$.

\begin{thm} Let $\alpha,\beta\in\cP$ be given. The  bijection in
Theorem \ref{main result} when restricted to
$\left(\M_{\A,\mathbb{N}}\times
\M_{\B,\mathbb{N}^\vee}\right)_{(\alpha,\beta)}$ gives a bijection
\begin{equation*}
\bigsqcup_{\lambda}SST_\A(\lambda/\alpha)\times
SST_\B(\lambda/\beta) \longrightarrow
\bigsqcup_{\eta}SST_\A(\beta/\eta)\times SST_\B(\alpha/\eta) \times
\M_{\A,\B}.
\end{equation*}
\end{thm}
\pf Since the bijection in Theorem \ref{main result} preserves the
plactic relations, we have a bijection
\begin{equation}\label{LRword-5}
\left(\M_{\A,\mathbb{N}}\times
\M_{\B,\mathbb{N}^\vee}\right)_{(\alpha,\beta)}{\longrightarrow}
\left( \M_{\B,\mathbb{N}^\vee}\times \M_{\A,\mathbb{N}}
\right)_{(\alpha,\beta)}\times \M_{\A,\B}.
\end{equation}
On the other hand, we have {\allowdisplaybreaks
\begin{align*}
&\left(\M_{\A,\mathbb{N}}\times
\M_{\B,\mathbb{N}^\vee}\right)_{(\alpha,\beta)} \\
&\stackrel{1-1}{\longleftrightarrow}\ \bigsqcup_{\mu, \nu}\
SST_\A(\mu)\times  SST_\B(\nu)\times \left(\B_{\mu}
\times \B_{\nu}^\vee\right)_{(\alpha,\beta)}\\
&\stackrel{1-1}{\longleftrightarrow} \bigsqcup_{\lambda,\mu, \nu}\
SST_\A(\mu)\times  SST_\B(\nu)\times {\bf LR}^{\lambda}_{\alpha
\mu}\times{\bf LR}^{\lambda}_{\beta \nu} \ \ \  \ \ \  \ \ \ \
\text{by (\ref{LRword-3})}
\\
&\stackrel{1-1}{\longleftrightarrow} \bigsqcup_{\lambda,\mu, \nu}\
SST_\A(\mu)\times {\bf LR}^{\lambda}_{ \mu \alpha}\times
SST_\B(\nu)\times\times{\bf LR}^{\lambda}_{ \nu\beta}\\
&\stackrel{1-1}{\longleftrightarrow}\
\bigsqcup_{\lambda}SST_\A(\lambda/\alpha)\times
SST_\B(\lambda/\beta)  \ \ \ \ \ \ \ \ \ \ \ \ \ \ \ \ \ \  \ \ \  \
\ \  \ \ \  \ \text{by (\ref{skewLR})}.
\end{align*}}
Similarly, we have
\begin{align*}
&\left( \M_{\B,\mathbb{N}^\vee}\times \M_{\A,\mathbb{N}}
\right)_{(\alpha,\beta)} \\
&\stackrel{1-1}{\longleftrightarrow}\ \bigsqcup_{\sigma, \tau}\
SST_\A(\sigma)\times  SST_\B(\tau)\times \left(\B_{\tau}^\vee
\times \B_{\sigma}\right)_{(\alpha,\beta)}\\
&\stackrel{1-1}{\longleftrightarrow} \bigsqcup_{\eta,\sigma, \tau}\
SST_\A(\sigma)\times  SST_\B(\tau)\times {\bf LR}^{\alpha}_{\eta
\tau}\times{\bf LR}^{\beta}_{\eta \sigma} \ \ \  \ \ \ \ \ \ \
\text{by (\ref{LRword-4})}
\\
&\stackrel{1-1}{\longleftrightarrow} \bigsqcup_{\eta,\sigma, \tau}\
SST_\A(\sigma)\times{\bf LR}^{\beta}_{ \sigma\eta} \times
SST_\B(\tau)\times {\bf
LR}^{\alpha}_{ \tau\eta}\\
&\stackrel{1-1}{\longleftrightarrow}\
\bigsqcup_{\eta}SST_\A(\beta/\eta)\times SST_\B(\alpha/\eta) \ \ \ \
\ \ \ \ \ \ \ \ \ \ \ \ \ \  \ \ \  \ \ \  \ \ \  \ \text{by
(\ref{skewLR})}.
\end{align*}}
Combining with (\ref{LRword-5}), we obtain the result. \qed

{\small

\end{document}